%% file: GPZ2023.tex
\documentclass[graybox]{svmult}

\usepackage{type1cm}        
\usepackage{makeidx}         
\usepackage{graphicx}        
\usepackage{multicol}        
\usepackage[bottom]{footmisc}

\usepackage{newtxtext}       %
\usepackage[varvw]{newtxmath}       

\makeindex             
\usepackage{tikz}
\usetikzlibrary{arrows}

\newcommand*\diff{\mathop{}\!\mathrm{d}}
\usepackage[normalem]{ulem}

\begin{document}

\title*{Asymptotic behavior of the heat semigroup on certain Riemannian
manifolds}

\author{Alexander Grigor'yan, Effie Papageorgiou and Hong-Wei Zhang}

\institute{
Alexander Grigor'yan \at Fakult\"{a}t f\"{u}r Mathematik, Universit\"{a}t Bielefeld, Postfach 100131, 33501 Bielefeld, Germany \at
\email{grigor@math.uni-bielefeld.de}
\and Effie Papageorgiou \at Department of Mathematics and Applied Mathematics, University of Crete, Crete, Greece \at
\email{papageoeffie@gmail.com}
\and Hong-Wei Zhang \at Department of Mathematics: Analysis, Logic and Discrete Mathematics Ghent University, Belgium \at
\email{hongwei.zhang@ugent.be}
}

\maketitle

\abstract{We show that, on a complete, connected and non-compact Riemannian 
manifold of non-negative Ricci curvature, the solution to the heat equation 
with $L^{1}$ initial data behaves asymptotically as the mass times the heat kernel. 
In contrast to the previously known results in negatively curved contexts, 
the radiality assumption on the initial data is not required. 
Similar long-time convergence results remain valid on more 
general manifolds satisfying the Li-Yau two-sided estimate of the heat kernel.
Moreover, we provide a counterexample such that this asymptotic phenomenon fails 
in sup norm on manifolds with two Euclidean ends.}


\section{Introduction}

\label{Section.1 Intro}

Let $\mathcal{M}$ be a Riemannian manifold of dimension $n\geq 2$ and $%
\Delta $ be the Laplace-Beltrami operator on $\mathcal{M}$. It is well
understood that the long-time behavior of solutions to the heat equation 
\begin{equation}
\partial _{t}u(t,x)\,=\,\Delta _{x}u(t,x),\qquad
\,u(0,x)\,=\,u_{0}(x),\qquad \,t>0,\,\,x\in \mathcal{M}  \label{S1 heat}
\end{equation}%
is very much related to the global geometry of $\mathcal{M}.$ This applies
also to the heat kernel $h_{t}\left( x,y\right) $ that is the minimal
positive fundamental solution of the heat equation or, equivalently, the
integral kernel of the heat semigroup $\exp \left( t\Delta \right) $ (see
for instance \cite{Dav1990,Sal2002,Gri2009}).

If the initial function $u_{0}$ belongs to the space $L^{p}\left(\mathcal{M},\mu
\right) $ with $p\in \lbrack 1,\infty )$ (where $\mu $ is the Riemannian
measure on $\mathcal{M}$) then the Cauchy problem (\ref{S1 heat}) has a
unique solution $u$ such that $u\left( t,\cdot \right) \in L^{p}$ for any $%
t>0$, and this solution is given by 
\begin{equation}
u\left( t,x\right) =\int_{\mathcal{M}}h_{t}\left( x,y\right) u_{0}\left( y\right) 
\mathrm{d}\mu (y).  \label{u}
\end{equation}%
The same is true for the case $p=\infty $ provided $\mathcal{M}$ is
stochastically complete. Hence, by a solution of (\ref{S1 heat}) we always
mean the function (\ref{u}).

The aim of this paper is to investigate the connection between the long-time
behavior of the solution $u\left( t,x\right) $ of (\ref{S1 heat}) and that
of the heat kernel $h_{t}\left( x,y\right) $. Let the initial function $%
u_{0} $ belong to $L^{1}(\mathcal{M}),$ and denote by $M=\int_{\mathcal{M}}u_{0}(x)\diff{\mu(x)}$ its mass. In the case when $\mathcal{M}=\mathbb{R}^{n}$
with the Euclidean metric, the heat kernel is given by 
\begin{equation*}
h_{t}(x,y)\,=\,(4\pi {t})^{-\frac{n}{2}}e^{-\frac{|x-y|^{2}}{4t}}
\end{equation*}%
and the solution to (\ref{S1 heat}) satisfies as $t\rightarrow \infty $ 
\begin{equation}
\Vert u(t,\,.\,)\,-\,Mh_{t}(\,.\,,x_{0})\Vert _{L^{1}(\mathbb{R}%
^{n})}\,\longrightarrow \,0  \label{S1 L1 R}
\end{equation}%
and 
\begin{equation}
t^{\frac{n}{2}}\,\Vert u(t,\,.\,)\,-\,Mh_{t}(\,.\,,x_{0})\Vert _{L^{\infty }(%
\mathbb{R}^{n})}\,\longrightarrow \,0.  \label{S1 Linf R}
\end{equation}%
By interpolation, a similar convergence holds with respect to any $L^{p}$
norm when $1<{p}<\infty $: 
\begin{equation*}
t^{\frac{n}{2p^{\prime }}}\,\Vert u(t,\,.\,)\,-\,Mh_{t}(\,.\,,0)\Vert
_{L^{p}(\mathbb{R}^{n})}\,\longrightarrow \,0
\end{equation*}%
where $p^{\prime }$ is the H\"{o}lder conjugate of $p$.

Note that (\ref{S1 L1 R}) and (\ref{S1 Linf R}) hold for \emph{any} choice
of $x_{0}$, which means that in the long run the solution $u\left( t,x\right) $
and the heat kernel $h_{t}\left( x,x_{0}\right) $ \textquotedblleft
forget\textquotedblright\ about the initial function $u_{0}$ resp. initial
point $x_{0}.$ We refer to a recent survey \cite{Vaz2018} for more details
about this property in the Euclidean setting.

The convergence properties (\ref{S1 L1 R}) and (\ref{S1 Linf R}) have an
interesting probabilistic meaning. Let $\left\{ X_{t}\right\} $ be Brownian
motion on $\mathcal{M}$ whose transition density is $h_{t}\left( x,y\right)
. $ Then (\ref{S1 Linf R}) means, in particular, that $X_{t}$ eventually
\textquotedblleft forgets\textquotedblright\ about its starting point $%
x_{0}, $ which corresponds to the fact that $X_{t}$ escapes to $\infty $
rotating chaotically in angular direction.

The situation is drastically different in hyperbolic spaces. It was shown by
V\'{a}zquez \cite{Vaz2019} that (\ref{S1 L1 R}) fails for a general initial
function $u_{0}\in L^{1}\left( \mathbb{H}^{n}\right) $ but is still true if $%
u_{0}$ is spherically symmetric around $x_{0}.$ Similar results were
obtained in \cite{APZ2023} in a more general setting of symmetric spaces of
non-compact type by using tools of harmonic analysis. Note that these spaces
have nonpositive sectional curvature. Recall that in hyperbolic spaces
Brownian motion $X_{t}$ tends to escape to $\infty $ along geodesics, which
means that it \textquotedblleft remembers\textquotedblright\ at least the
direction of the starting point $x_{0}.$

Our main result is the following theorem that deals with manifolds of
non-negative Ricci curvature. Denote by $B\left( x,r\right) $ the geodesic
ball on $\mathcal{M}$ of radius $r$ centered at $x\in \mathcal{M}$ and set $%
V\left( x,r\right) =\mu \left( B\left( x,r\right) \right) .$

\begin{theorem}
\label{S1 Main thm} Let $\mathcal{M}$ be a complete, connected and
non-compact Riemannian manifold of non-negative Ricci curvature. Fix a base
point $x_{0}\in \mathcal{M}$ and suppose that $u_{0}\in {L}^{1}(\mathcal{M})$%
. Then the solution to the heat equation (\ref{S1 heat}) satisfies as $%
t\rightarrow \infty $ 
\begin{equation}
\Vert u(t,\,.\,)\,-\,Mh_{t}(\,.\,,x_{0})\Vert _{L^{1}(\mathcal{M}%
)}\,\rightarrow \,0  \label{S1 Main thm L1}
\end{equation}%
and%
\begin{equation}
\,\Vert \left\vert u(t,\,.\,)\,-\,Mh_{t}(\,.\,,x_{0})\right\vert V(\cdot ,%
\sqrt{t})\Vert _{L^{\infty }(\mathcal{M})}\,\rightarrow 0.
\label{S1 Main thm Linf}
\end{equation}
\end{theorem}

\begin{remark}
By interpolation between (\ref{S1 Main thm L1}) and (\ref{S1 Main thm Linf}%
), we obtain for any $p\in \left( 1,\infty \right) $ 
\begin{equation}
\,\Vert \left\vert u(t,\,.\,)\,-\,Mh_{t}(\,.\,,x_{0})\right\vert V(\cdot ,%
\sqrt{t})^{1/p^{\prime }}\Vert _{L^{p}(\mathcal{M})}\,\longrightarrow \,0%
\mathnormal{.}  \label{S1 Main thm Lp}
\end{equation}
\end{remark}

In Section \ref{Section.2 Kernel} we will give a short review of 
different estimates of the heat kernel and reformulate
Theorem \ref{S1 Main thm} in more general ways.
In Section \ref{Sec3} we prove Theorem \ref{S1
Main thm}. An essential idea of the proof is to describe the critical region
where the heat kernel concentrates. 
In the last section, we show that an estimate as 
\eqref{S1 Main thm Linf} fails on manifolds with two Euclidean ends.

\section{Auxiliary results}\label{Section.2 Kernel}
\subsection{Heat kernel estimates}
From now on, $\mathcal{M}$ denotes a complete,
connected, non-compact Riemannian manifold of dimension $n\geq 2$. Let $\mu $
be the Riemannian measure on $\mathcal{M}$. Let $d(x,y)$ be the geodesic
distance between two points $x,y\in \mathcal{M}$, and $V(x,r)=\mu \left(
B\left( x,r\right) \right) $ be the Riemannian volume of the geodesic ball $%
B(x,r)$ of radius $r$ centered at $x\in \mathcal{M}$.

Throughout the paper we follow the convention that $C,C_{1},...$ denote
large positive constants whereas $c,c_{1},...$ are small positive constants.
These constants may depend on $\mathcal{M}$ but do not depend on the
variables $x,y,t$. Moreover, the notation $A\lesssim {B}$ between two
positive expressions means that $A\leq {C}B$, and $A\asymp {B}$ means $%
cB\leq {A}\leq {C}B$.

We say that $\mathcal{M}$ satisfies the volume doubling property if, for all 
$x\in \mathcal{M}$ and $r>0$, we have 
\begin{equation}
V(x,2r)\,\leq \,C\,V(x,r).  \label{VD}
\end{equation}%
It follows from (\ref{VD}) that there exist some positive constants $\nu
,\nu ^{\prime }>0$ such that 
\begin{equation}
c\,\left( \frac{R}{r}\right) ^{\nu ^{\prime }}\,\leq \,\frac{V(x,R)}{V(x,r)}%
\,\leq \,C\,\left( \frac{R}{r}\right) ^{\nu }  \label{S2 Comp same center}
\end{equation}%
for all $x\in \mathcal{M}$ and $0<r\leq {R}$ (see for instance \cite[Section
15.6]{Gri2009}). Moreover, (\ref{S2 Comp same center}) implies that, for all 
$x,y\in \mathcal{M}$ and $r>0,$ 
\begin{equation}
\frac{V\left( x,r\right) }{V\left( y,r\right) }\leq C\left( 1+\frac{d\left(
x,y\right) }{r}\right) ^{\nu }.  \label{S2 Comp diff center}
\end{equation}

The integral kernel $h_{t}\left( x,y\right) $ of the heat semigroup $\exp
(t\Delta )$ is the smallest positive fundamental solution to the heat
equation (\ref{S1 heat}). It is known that $h_{t}\left( x,y\right) $ is
smooth in $\left( t,x,y\right) $, symmetric in $x,y$, and
satisfies the semigroup identity (see for instance \cite{Str1983}, \cite%
{Gri2009}). Besides, for all $y\in\mathcal{M}$ and $t>0$ 
\begin{equation*}
\int_{\mathcal{M}}h_{t}\left( x,y\right) \mathrm{d}\mu (x)\leq 1.
\end{equation*}%
The manifold $\mathcal{M}$ is called stochastically complete if for all $%
y\in\mathcal{M}$ and $t>0$ 
\begin{equation*}
\int_{\mathcal{M}}h_{t}(x,y)\,\mathrm{d}\mu (x)=\,1.
\end{equation*}%
It is known that if $\mathcal{M}$ is geodesically complete and, for some $%
x_{0}\in \mathcal{M}$ and all large enough $r,$%
\begin{equation*}
V\left( x_{0},r\right) \leq e^{Cr^{2}},
\end{equation*}%
then $\mathcal{M}$ is stochastically complete. In particular, the volume
doubling property (\ref{VD}) implies that $\mathcal{M}$ is stochastically
complete.

When the Ricci curvature of $\mathcal{M}$ is non-negative, 
the following two-sided estimates of the heat kernel were proved by Li and
Yau \cite{LiYa1986}:
\begin{equation}
\frac{c_{1}}{V(x,\sqrt{t})}\,\exp \Big(-C_{1}\,\frac{d^{2}(x,y)}{t}\Big)%
\,\leq \,h_{t}(x,y)\,\leq \,\frac{C_{2}}{V(x,\sqrt{t})}\,\exp \Big(-c_{2}\,%
\frac{d^{2}(x,y)}{t}\Big). \label{S2 twosided}
\end{equation}%
Besides, Li and Yau proved in \cite{LiYa1986} also the following gradient
estimate for any positive solution $u\left( t,x\right) $ of the heat
equation $\partial _{t}u=\Delta u$ on $\mathbb{R}_{+}\times \mathcal{M}$: 
\begin{equation}
\frac{|\nabla u|^{2}}{u^{2}}\,-\,\frac{\partial _{t}u}{u}\,\leq \,\frac{C}{t}
\label{S2 LY gradient}
\end{equation}%
with $C=\frac{n}{2}.$ By a result of \cite{Gri1995}, the upper bound of the
heat kernel in (\ref{S2 twosided}), that is, 
\begin{equation}
h_{t}(x,y)\,\leq \,\frac{C}{V(x,\sqrt{t})}\,\exp \Big(-c\,\frac{d^{2}(x,y)}{t%
}\Big),  \label{S2 upper}
\end{equation}%
implies the following estimate of the time derivative:%
\begin{equation}
\Big|\frac{\partial {h_{t}}}{\partial {t}}(x,y)\Big|\,\leq \,\frac{C}{t\,V(x,%
\sqrt{t})}\,\exp \Big(-c\frac{d^{2}(x,y)}{t}\Big).
\label{S2 time derivative}
\end{equation}%
It follows from (\ref{S2 LY gradient}) that%
\begin{equation*}
\left\vert \nabla u\right\vert ^{2}\leq u\partial _{t}u+\frac{C}{t}u^{2}
\end{equation*}%
and, hence, 
\begin{equation*}
\left\vert \nabla u\right\vert \leq \sqrt{u\left\vert \partial
_{t}u\right\vert }+\frac{C}{\sqrt{t}}u.
\end{equation*}%
Applying this for the function $u\left( t,y\right) =h_{t}\left( x,y\right) $
and combining with (\ref{S2 upper}) and (\ref{S2 time derivative}), we
obtain that 
\begin{equation}
|\nabla _{y}\,h_{t}(x,y)|\,\leq \,\frac{C}{\sqrt{t}\,V(x,\sqrt{t})}\,\exp %
\Big(-c\,\frac{d^{2}(x,y)}{t}\Big).  \label{S2 gradient}
\end{equation}%

It is known that the upper bound \eqref{S2 upper} of $h_t$
(and, hence, its consequence \eqref{S2 time derivative}) holds on a larger class 
of Riemannian manifolds satisfying a so called \textit{relative Faber-Krahn inequality}
(see, for example, \cite{Gri2009}). 
But the estimate \eqref{S2 gradient} of the gradient $\nabla{h_t}$
is much more subtle and requires more serious hypotheses, for example,
non-negative Ricci curvature as we consider here. 

Let us observe that the gradient estimate \eqref{S2 gradient}
is a particular (and limiting) case of the following Hölder estimate:
\begin{align}
    |h_{t}(x,y)\,-\,h_{t}(x,z)|\,
    \le\,
    \Big(\frac{d(y,z)}{\sqrt{t}}\Big)^{\theta}\,
    \frac{C}{V(x,\sqrt{t})}\,
    \exp\Big(-c\frac{d^2(x,y)}{t}\Big),
    \label{Holder}
\end{align}
for some $0<\theta\le1$, and all $t>0$, $x,y,z\in\mathcal{M}$ 
such that $d(y,z)\le\sqrt{t}$. It is important to mention that
\eqref{Holder} is a consequence of the two-sided estimate \eqref{S2 twosided}
alone (for some $0<\theta<1$, see \cite[Theorem 5.4.12]{Sal2002}), and, hence,
\eqref{Holder} holds on more general classes of manifolds
than \eqref{S2 gradient}.

Finally, let us observe that on spaces of essentially negative curvature the
above estimates of the heat kernel typically fail: for example, these are
hyperbolic spaces \cite{DaMa1988}, non-compact symmetric spaces \cite{AnOs2003}, 
asymptotically hyperbolic manifolds \cite{ChHa2020}, and fractal-like manifolds 
\cite{Bar1998}.

\subsection{Alternative statements of the main theorem}

Now let us reformulate Theorem \ref{S1 Main thm} in a bit more general way.

\begin{theorem}
\label{Tmain}Let $\mathcal{M}$ be a geodesically complete non-compact
manifold that satisfies the following conditions:

\begin{itemize}
\item the volume doubling condition (\ref{VD});

\item the upper bound$\,$(\ref{S2 upper}) of the heat kernel;

\item the Hölder regularity \eqref{Holder} of the heat kernel.
\end{itemize}

Then the conclusions of Theorem \ref{S1 Main thm} hold.
\end{theorem}

Apart from manifolds with non-negative Ricci curvature, the
above-described manifolds cover many other examples.
Let us recall that, on a complete Riemannian manifold,
the following three properties are equivalent:

\begin{itemize}
\item The two-sided estimate \eqref{S2 twosided} of the heat kernel;

\item
The uniform parabolic Harnack inequality:
\begin{align}
    \sup_{(\frac{T}{4},\frac{T}{2})\times{B(x,\frac{r}{2})}}u(t,x)\;
    \le\,C\,
    \sup_{(\frac{3T}{4},T)\times{B(x,\frac{r}{2})}}u(t,x),
    \label{PHI}
\end{align}
where $u(t,x)$ is a non-negative solution of the heat equation 
$\partial _{t}u=\Delta u$ in a cylinder $(0,T)\times{B(x,r)}$ 
with $x\in\mathcal{M}$, $r>0$ and $T=r^{2}$.

\item 
The conjunction of the volume doubling property \eqref{VD} and the
Poincaré inequality:
\begin{align}
    \int_{B(x,r)}|f-f_{B}|^{2}\,\diff{\mu}\,
    \le\,C\,r^{2}\,
     \int_{B(x,r)}|\nabla{f}|^{2}\,\diff{\mu},
     \label{PI}
\end{align}
for all $x\in\mathcal{M}$, $t>0$, and bounded Lipschitz functions $f$ in
$B(x,r)$. Here, $f_{B}$ is the mean of $f$ over $B(x,r)$.
\end{itemize}
See, for instance, \cite {FaSt1986,Gri1991,Sal1992,Sal2002} for more details.
Manifolds satisfying these equivalent conditions include complete manifolds 
with non-negative Ricci curvature, connected Lie groups with polynomial 
volume growth, co-compact covering manifolds whose deck transformation 
has polynomial growth, and many others. 
We refer to \cite[pp. 417--418]{Sal2010} for a list of examples. 
These equivalent conditions yield particularly the Hölder regularity 
\eqref{Holder} for some $0<\theta<1$, 
see \cite[Theorem 5.4.12]{Sal2002}. Hence the following
corollary follows.

\begin{corollary}
Let $\mathcal{M}$ be a geodesically complete non-compact
manifold that satisfies one of the following equivalent conditions:
\begin{itemize}
\item the two-sided estimate \eqref{S2 twosided} of the heat kernel;

\item the uniform parabolic Harnack inequality \eqref{PHI};

\item the conjunction of the volume doubling property \eqref{VD} and the
Poincaré inequality \eqref{PI}.
\end{itemize}

Then the conclusions of Theorem \ref{S1 Main thm} hold.
\end{corollary}

In fact, we do not need a differential operator like the Laplacian (or
sub-Laplacian). All we need is a function $h_t(x,y)$ satisfying certain
estimates, and the result can be formulated as a property of integral
operators. For example, our method can be applied on a metric space $(X,d,\mu)$ 
of homogeneous type, assuming that one has a semigroup of self-adjoint 
linear operators acting on $L^{2}(X,\mu)$ and admitting a good integral kernel, 
see for instance \cite{DzPr2018}.

\subsection{Two preliminary lemmas}

We prove Theorem \ref{Tmain} in the next section. In the next two lemmas we
describe some consequences of the hypotheses of Theorem \ref{Tmain}, in
particular, the critical annulus where the heat kernel concentrates.

\begin{lemma}
\label{S3 Lemma1}Under the volume doubling condition (\ref{VD}), we have, 
for any $c>0$, $x_{0}\in\mathcal{M}$, and $t>0,$%
\begin{equation}
\int_{\mathcal{M}}\,\exp \left( -\,c\,\frac{%
d\left( x,x_{0}\right) ^{2}}{t}\right)\, \frac{\diff{\mu} (x)}{V(x_{0},\sqrt{t})}\,\lesssim 1  \label{int}
\end{equation}%
and, for any $N\in \mathbb{N}$ and all $r\geq \sqrt{t}$,%
\begin{equation}
\int_{B\left( x_{0},r\right) ^{c}}
\,\exp \left( -\,c\,\frac{d\left( x,x_{0}\right) ^{2}}{t}\right) \,
\frac{\diff{\mu} (x)}{V(x_{0},\sqrt{t})}\,\lesssim
\left( \frac{r}{\sqrt{t}}\right) ^{-N}.  \label{intr}
\end{equation}
\end{lemma}

\begin{proof}
Let us prove first (\ref{intr}). Using (\ref{S2 Comp same center}),
we have
\begin{align*}
&\frac{1}{V(x_{0},\sqrt{t})}\,
\int_{d(x,x_0)\ge{r}}
\,\exp \left( -\,c\,\frac{d\left( x,x_{0}\right) ^{2}}{t}\right) \,
\diff{\mu(x)} \\[5pt]
& =\,\frac{1}{V(x_{0},\sqrt{t})}\,\sum_{j=0}^{\infty
}\,\int_{2^{j}r\,\leq d(x,x_{0})\,<\,2^{j+1}r}\,
\exp \left( -\,c\,\frac{d\left( x,x_{0}\right) ^{2}}{t}\right)\,
\diff{\mu} (x) \\%
[5pt]
& \lesssim \,\sum_{j=0}^{\infty }\,\frac{V(x_{0},2^{j+1}r)}{V(x_{0},\sqrt{t})%
}\,\exp \left( -\,c\,\frac{\left( 2^{j}r\right) ^{2}}{t}\right) \\[5pt]
& \lesssim \,\sum_{j=0}^{\infty }\,\left( \frac{2^{j+1}r}{\sqrt{t}}\right)
^{\nu }\,\exp \left( -\,c\,\frac{2^{2j}r^{2}}{t}\right) .
\end{align*}%
Since $\exp \left( -s\right) \leq \left( N!\right) s^{-N}$ for any $s>0$ and 
$N\in \mathbb{N}$, we obtain that 
\begin{align*}
\int_{B\left( x_{0},r\right) ^{c}}
\,\exp \left( -\,c\,\frac{d\left( x,x_{0}\right) ^{2}}{t}\right) \,
\frac{\diff{\mu} (x)}{V(x_{0},\sqrt{t})}\,
&\lesssim
\,\sum_{j=0}^{\infty }\,\left( \frac{r}{\sqrt{t}}\right) ^{\nu }2^{\nu
j-2Nj}\left( \frac{t}{r^{2}}\right) ^{N}\\[5pt]
&\lesssim \left( \frac{r}{\sqrt{t}}%
\right) ^{\nu -2N},
\end{align*}%
which proves (\ref{intr}).

In order to prove (\ref{int}) we apply (\ref{intr}) with $r=\sqrt{t}$ and $%
N=1$ and obtain%
\begin{align*}
&\int_{\mathcal{M}}\,\exp \left( -\,c\,\frac{%
d\left( x,x_{0}\right) ^{2}}{t}\right)\,
\frac{\diff{\mu} (x)}{V(x_{0},\sqrt{t})}\\[5pt]
&=\left( \int_{B\left(
x_{0},\sqrt{t}\right) ^{c}}+\int_{B\left( x_{0},\sqrt{t}\right) }\right)\,
\exp \left( -\,c\,\frac{d\left( x,x_{0}\right) ^{2}%
}{t}\right)\,
\frac{\diff{\mu} (x)}{V(x_{0},\sqrt{t})}\\[5pt]
&\lesssim 1+\int_{B\left( x_{0},\sqrt{t}\right) }\frac{\diff{\mu} (x)}{%
V(x_{0},\sqrt{t})}=2.
\end{align*}
\hspace{1pt}
\end{proof}

The next lemma describes the annulus where the heat kernel concentrates. Let
us fix a point $x_{0}\in \mathcal{M}$ and a positive function $\varphi
\left( t\right) $ such that $\varphi \left( t\right) \rightarrow 0$ and $\varphi(t)\sqrt{t}\rightarrow\infty$ as $%
t\rightarrow \infty $. For any $t>0$, define the following annulus in $%
\mathcal{M}$:%
\begin{equation}
\Omega _{t}\,=\left\{ \,{x\in \mathcal{M}\,\ \big|\ \,\varphi (t)\sqrt{t}%
\,\leq \,d(x,x_{0})\,\leq \,}\frac{\sqrt{t}}{\varphi (t)}\right\} .
\label{Om}
\end{equation}

\begin{lemma}
\label{S2 Lemma concentration} Under the hypotheses (\ref{VD}) and $\,$(\ref%
{S2 upper}), we have for all large enough $t$ 
\begin{equation}
\int_{\mathcal{M}\setminus \Omega _{t}}\,h_{t}(x,x_{0})\diff{\mu(x)}%
\,\lesssim \,\varphi \left( t\right) ^{\nu ^{\prime }}  \label{M-Om}
\end{equation}%
where $\nu ^{\prime }$ is the exponent from (\ref{S2 Comp same center}).
Consequently, 
\begin{equation}
\int_{\Omega _{t}}\,h_{t}(x,x_{0})\diff{\mu(x)}\,\longrightarrow \,1\ \ \ 
\text{\textnormal{as}}\,\,\,t\longrightarrow \infty .  \label{1}
\end{equation}
\end{lemma}

\begin{proof}
Since $\mathcal{M}$ is stochastically complete, (\ref{1}) follows from (\ref%
{M-Om}) and $\varphi \left( t\right) \rightarrow 0.$ Since 
\begin{equation*}
\mathcal{M}\setminus \Omega _{t}=B\left( x_{0},\varphi(t)\sqrt{t}\right)
\cup B^{c}\left( x_{0},\sqrt{t}/\varphi(t)\right),
\end{equation*}%
we estimate the integrals over $B( x_{0},\varphi (t)\sqrt{t}) $
and $B^{c}( x_{0},\sqrt{t}/\varphi(t)) $ separately. Assume that $t$
is large enough so that $\varphi \left( t\right) <1.$ Using (\ref{S2 upper})
and (\ref{S2 Comp same center}) we obtain 
\begin{equation*}
\int_{d(x,x_{0})<\varphi (t)\sqrt{t}}\,h_{t}(x,x_{0})\diff{\mu(x)}\,\lesssim
\,\frac{V(x_{0},\varphi (t)\sqrt{t})}{V(x_{0},\sqrt{t})}\,\lesssim \,\varphi
(t)^{\nu ^{\prime }}.
\end{equation*}%
Using (\ref{S2 upper}) and (\ref{intr}) with $r=\frac{\sqrt{t}}{\varphi (t)}$
and any $N\in \mathbb{N}$, we obtain%
\begin{eqnarray*}
\int_{d(x,x_{0})>\frac{\sqrt{t}}{\varphi (t)}}\,h_{t}(x,x_{0})\diff{\mu}%
(x)\, &\lesssim &\int_{B\left( x_{0},r\right) ^{c}}\frac{\diff{\mu} (x)}{%
V(x_{0},\sqrt{t})}\,\exp \left( -\,c\,\frac{d\left( x,x_{0}\right) ^{2}}{t}%
\right) \\
&\lesssim &\left( \frac{r}{\sqrt{t}}\right) ^{-N}=\varphi \left( t\right)
^{N},
\end{eqnarray*}%
whence the claim follows.
\end{proof}

\section{Proof of the main theorem}

\label{Sec3}We start the proof of Theorem \ref{Tmain} with continuous
compactly supported initial data and prove the asymptotic properties of the
solution in the $L^{1}$ norm, working separately outside and inside the
critical region $\Omega _{t}$. Then, we show that these properties remain
valid for all $L^{1}$ initial data by using a density argument. The $%
L^{\infty }$ convergence is proved in the same spirit.

\begin{proposition}
\label{S3 proposition}
Let $\mathcal{M}$ be a manifold defined as in Theorem \ref{Tmain}.
Let $x_{0}\in \mathcal{M}$ and $u_{0}\in \mathcal{C}%
_{c}(B(x_{0},a))$ for some $a>0$. Assume that $\varphi \left( t\right) $ is
a positive function such that $\varphi \left( t\right) \rightarrow 0$ and $%
\varphi \left( t\right) \sqrt{t}\rightarrow \infty $ as $t\rightarrow \infty
.$ Then the solution (\ref{u}) satisfies 
\begin{equation}
\Vert u(t,\,.\,)\,-\,Mh_{t}(\,.\,,x_{0})\Vert _{L^{1}(\mathcal{M}%
\smallsetminus \Omega _{t})}\,\lesssim \,\varphi (t)^{\nu ^{\prime }}
\label{S3 prop outside}
\end{equation}%
and 
\begin{equation}
\Vert u(t,\,.\,)\,-\,Mh_{t}(\,.\,,x_{0})\Vert _{L^{1}(\Omega
_{t})}\,\lesssim \,t^{-\frac{\theta}{2}},  \label{S3 prop inside}
\end{equation}%
where $\Omega _{t}$ is defined by (\ref{Om}), $\nu ^{\prime }$ and $\theta$
are exponents from (\ref{S2 Comp same center}) and \eqref{Holder}, and 
\begin{equation*}
M=\int_{\mathcal{M}}\,u_{0}(x)\diff{\mu(x)}.
\end{equation*}%
Consequently, 
\begin{equation}
\Vert u(t,\,.\,)\,-\,Mh_{t}(\,.\,,x_{0})\Vert _{L^{1}\left( \mathcal{M}%
\right) }\lesssim t^{-\eta }  \label{eta}
\end{equation}%
for any $\eta <\min (\nu ^{\prime },\theta)/2.$
\end{proposition}

\begin{proof}
By the upper bound$\,$(\ref{S2 upper}) of the heat kernel and
Lemma \ref{S2 Lemma concentration} we have 
\begin{equation*}
\Vert h_{t}(\,.\,,x_{0})\Vert _{L^{1}(\mathcal{M}\smallsetminus \Omega
_{t})}\,\lesssim \,\varphi (t)^{\nu ^{\prime }}
\end{equation*}%
so that (\ref{S3 prop outside}) will follow if we prove that 
\begin{equation}
\Vert u(t,\,.\,)\Vert _{L^{1}(\mathcal{M}\smallsetminus \Omega
_{t})}\,\lesssim \,\varphi (t)^{\nu ^{\prime }}.  \label{S3 u outside}
\end{equation}%
We write 
\begin{align*}
\Vert u(t,\,.\,)\Vert _{L^{1}(\mathcal{M}\smallsetminus \Omega _{t})}\,&
=\,\int_{\mathcal{M}\smallsetminus \Omega _{t}}\,\left\vert
\int_{B(x_{0},a)}\,h_{t}(x,y)\,u_{0}(y)\diff{\mu(y)}\right\vert\diff{\mu(x)}
\\[5pt]
& \leq \,\int_{B(x_{0},a)}\,|u_{0}(y)|\,\Big\lbrace{\int_{\mathcal{M}%
\smallsetminus \Omega _{t}}\,h_{t}(x,y)\diff{\mu(x)}}\Big\rbrace\,%
\diff{\mu(y)}.
\end{align*}%
Notice that $x\in \mathcal{M}\smallsetminus \Omega _{t}$ and $y\in {%
B(x_{0},a)}$ imply $x\in \mathcal{M}\smallsetminus \widetilde{\Omega }_{t,y}$%
, where 
\begin{equation*}
\widetilde{\Omega }_{t,y}\,=\,\left\{ {x\in \mathcal{M}\,\big|\,2\,\varphi
(t)\sqrt{t}\,\leq \,d(x,y)\,\leq \,\frac{1}{2}\,\frac{\sqrt{t}}{\varphi (t)}}%
\right\}
\end{equation*}%
provided $t$ is large enough. Indeed, if $x\in \widetilde{\Omega }_{t,y}$
then 
\begin{equation*}
d(x,x_{0})\,\leq \,d(x,y)\,+\,d(y,x_{0})\,\leq \,\frac{1}{2}\,\frac{\sqrt{t}%
}{\varphi (t)}\,+\,a\,\leq \,\frac{\sqrt{t}}{\varphi (t)}
\end{equation*}%
and 
\begin{equation*}
d(x,x_{0})\,\geq \,d(x,y)\,-\,d(y,x_{0})\,\geq \,{2\,\varphi (t)\sqrt{t}\,}%
-\,a\geq \,\,\varphi (t)\sqrt{t}\,
\end{equation*}%
for $t$ large enough, since $\varphi (t)\rightarrow 0$ and $\varphi (t)\sqrt{%
t}\rightarrow \infty $ as $t\rightarrow \infty $. It follows that $x\in
\Omega _{t}$.

Applying (\ref{M-Om}) with $\widetilde{\Omega }_{t,y}$ instead of $\Omega
_{t}$ we obtain 
\begin{equation*}
\int_{\mathcal{M}\smallsetminus \Omega _{t}}\,h_{t}(x,y)\,\diff{\mu(x)}%
\,\leq \,\int_{\mathcal{M}\smallsetminus \widetilde{\Omega }%
_{t,y}}\,h_{t}(x,y)\diff{\mu(x)}\,\lesssim \,\varphi (t)^{\nu ^{\prime }},
\end{equation*}%
whence (\ref{S3 u outside}) follows.

Now, let us turn to (\ref{S3 prop inside}). Observe that 
\begin{align}
u(t,x)\,-\,Mh_{t}(x,x_{0})\,& =\,\int_{\mathcal{M}}\,u_{0}(y)\,\left(
h_{t}(x,y)\,-\,h_{t}(x,x_{0})\right)\diff{\mu(y)}  \notag \\
& =\int_{B\left( x_{0},a\right) }\,u_{0}(y)\,\left(
h_{t}(x,y)\,-\,h_{t}(x,x_{0})\right)\diff{\mu(y)} .  \label{comment}
\end{align}%

One the one hand, we deduce from the Hölder regularity \eqref{Holder} that
\begin{align}
&|u(t,x)\,-\,Mh_{t}(x,x_{0})|\notag\\[5pt]
&\le\,
\frac{C}{V(x,\sqrt{t})}\,\exp\Big( -c\,\frac{d^{2}(x,x_{0})}{t}\Big)\,
\int_{B(x_{0},a)}\,\Big(\frac{d(x_0,y)}{\sqrt{t}}\Big)^{\theta}\,
\left\vert u_{0}(y)\right\vert\diff{\mu(y)}\notag\\[5pt]
&\lesssim\,
t^{-\frac{\theta}{2}}\,
\frac{1}{V(x,\sqrt{t})}\,\exp\Big( -c\,\frac{d^{2}(x,x_{0})}{t}\Big),
\label{S3 difference}
\end{align}
for some $0<\theta\le1$ and for $t$ large enough such that 
$d(x_0,y)\le{a}\le\sqrt{t}$.
On the other hand, (\ref{S2 Comp diff center}) implies that
\begin{align}
\frac{1}{V(x,\sqrt{t})}\exp \left( -c\,\frac{d^{2}(x,x_{0})}{t}\right) \,&
\lesssim \,\,\left( \frac{d(x,x_{0})}{\sqrt{t}}+1\right) ^{\nu }\,\frac{1}{%
V(x_{0},\sqrt{t})}\exp \left( -c\,\frac{d^{2}(x,x_{0})}{t}\right)  \notag \\%
[0.07in]
&\lesssim \,\frac{1}{V(x_{0},\sqrt{t})}\,\exp \left( -c\,\frac{%
d^{2}(x,x_{0})}{2t}\right) .  \label{S3 proof ingredient 2}
\end{align}%

Substituting (\ref{S3 proof ingredient 2}) into the right-hand side of (\ref%
{S3 difference}) and integrating in $x$ over $\Omega _{t}$, we obtain by (%
\ref{int}) 
\begin{equation*}
\int_{\Omega _{t}}\,|u(t,x)\,-\,Mh_{t}(x,x_{0})|\diff{\mu(x)}\,\lesssim \,%
t^{-\frac{\theta}{2}}\,\int_{\Omega _{t}}\frac{\diff{\mu} (x)}{V(x_{0},\sqrt{t})%
}\exp \left( -c\,\frac{d^{2}(x,x_{0})}{2t}\right) 
\lesssim \,t^{-\frac{\theta}{2}}\,.
\end{equation*}

Finally, (\ref{eta}) follows by adding up (\ref{S3 prop outside}) and (\ref%
{S3 prop inside}) with $\varphi \left( t\right) =t^{\varepsilon -\frac{\theta}{2}%
} $ with small enough $\varepsilon .$
\end{proof}

Next, we prove our main theorem.\label{here}

\begin{proof}[Proof of Theorem \protect\ref{Tmain}]
Given $u_{0}\in {L^{1}(\mathcal{M})}$, fix $\varepsilon >0$ and choose $%
\widetilde{u}_{0}\in \mathcal{C}_{c}(\mathcal{M})$ such that 
\begin{equation*}
\Vert u_{0}-\widetilde{u}_{0}\Vert _{L^{1}(\mathcal{M})}<\frac{\varepsilon }{%
3}.
\end{equation*}%
Let us prove first (\ref{S1 Main thm L1}). Setting $\widetilde{M}=\int_{%
\mathcal{M}}\,\widetilde{u}_{0}(y)\diff{\mu(y)}$ we have 
\begin{equation}
|M\,-\,\widetilde{M}|\,\leq \,\int_{\mathcal{M}}\diff{\mu(y)}\,|u_{0}(y)\,-\,%
\widetilde{u}_{0}(y)|\,=\,\Vert u_{0}-\widetilde{u}_{0}\Vert _{L^{1}(%
\mathcal{M})}\,<\,\frac{\varepsilon }{3}.  \label{S3 mass diff}
\end{equation}%
It follows that 
\begin{equation}
\Vert Mh_{t}(\,.\,,x_{0})\,-\,\widetilde{M}h_{t}(\,.\,,x_{0})\Vert _{L^{1}(%
\mathcal{M})}\,\leq \,|M\,-\,\widetilde{M}|\,\Vert h_{t}(\,.\,,x_{0})\Vert
_{L^{1}(\mathcal{M})}<\,\frac{\varepsilon }{3}.  \label{S3 density1}
\end{equation}%
Let $\widetilde{u}(t,x)$ be the semigroup solution to the heat equation with
the initial data $\widetilde{u}_{0}$. Then we have 
\begin{align}
&\Vert u(t,\,.\,)\,-\,\widetilde{u}(t,\,.\,)\Vert _{L^{1}(\mathcal{M})}\notag\\[5pt]
&\leq\,\int_{\mathcal{M}}\,|u_{0}(y)\,-\,\widetilde{u}_{0}(y)|\,\Big\lbrace{\int_{\mathcal{M}}\diff{\mu(x)}\,h_{t}(x,y)}\Big\rbrace \diff{\mu(y)}\,
<\,\frac{\varepsilon }{3}.  \label{S3 density2}
\end{align}%
By the estimate (\ref{eta}) of Proposition \ref{S3 proposition}, we have,
for sufficiently large $t,$ 
\begin{equation}
\Vert \widetilde{u}(t,\,.\,)\,-\,\widetilde{M}h_{t}(\,.\,,x_{0})\Vert
_{L^{1}(\mathcal{M})}\,<\,\frac{\varepsilon }{3}.  \label{S3 density3}
\end{equation}%
Combining (\ref{S3 density1}), (\ref{S3 density2}) and (\ref{S3 density3})
together, we obtain 
\begin{equation*}
\Vert u(t,\,.\,)\,-\,Mh_{t}(\,.\,,x_{0})\Vert _{L^{1}(\mathcal{M}%
)}\,<\,\varepsilon
\end{equation*}%
which finishes the proof of (\ref{S1 Main thm L1}).

Let us turn to the sup norm convergence (\ref{S1 Main thm Linf}). Since $%
\widetilde{u}_{0}$ is compactly supported, we can apply (\ref{S3 difference}%
) which gives for all $x$%
\begin{equation}
|\widetilde{u}(t,x)\,-\,\widetilde{M}h_{t}(x,x_{0})|\,\leq \frac{C}{\sqrt{t}%
V(x,\sqrt{t})}\leq \frac{\varepsilon }{V(x,\sqrt{t})}
\label{S3 density inf1}
\end{equation}%
for sufficiently large $t.$ Next, we have 
\begin{align}
|u(t,x)\,-\,\widetilde{u}(t,x)|\,&\leq \,\int_{\mathcal{M}}\,h_{t}(x,y)\,%
\left\vert u_{0}(y)-\widetilde{u}_{0}(y)\right\vert\diff{\mu(y)}\notag\\[5pt]
&\leq\,
\frac{C\left\Vert u_{0}-\widetilde{u}_{0}\right\Vert _{L^{1}}}{V(x,\sqrt{t})}%
\lesssim \,\frac{\varepsilon }{V(x,\sqrt{t})}  \label{S3 density inf2}
\end{align}%
and 
\begin{equation}
|Mh_{t}(x,x_{0})\,-\,\widetilde{M}h_{t}(x,x_{0})|\,=\,|M-\widetilde{M}%
|\,h_{t}( x,x_{0}) \,\lesssim \frac{\varepsilon }{V( x,\sqrt{%
t}) }.  \label{S3 density inf3}
\end{equation}%
Putting together (\ref{S3 density inf1}), (\ref{S3 density inf2}) and (\ref%
{S3 density inf3}), we conclude that 
\begin{equation*}
\left\vert \,u(t,\,x\,)\,-\,Mh_{t}(\,x\,,x_{0})\right\vert \lesssim \,\frac{%
\varepsilon }{V(x,\sqrt{t}) }
\end{equation*}%
whence (\ref{S1 Main thm Linf}) follows.
\end{proof}


\section{A counterexample}

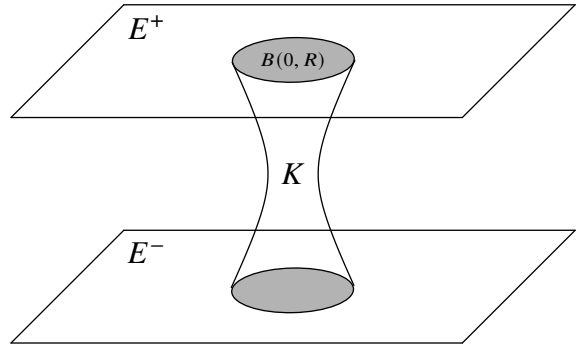
\begin{figure}
\begin{tikzpicture}[scale=1.5,line cap=round,line join=round,>=triangle 45,x=1cm,y=1cm]
\clip(-4.5,-1.25) rectangle (4,2.1);
\draw [line width=0.5pt] (-2,2)-- (2,2);
\draw [line width=0.5pt] (2,2)-- (1,1);
\draw [line width=0.5pt] (1,1)-- (-3,1);
\draw [line width=0.5pt] (-3,1)-- (-2,2);
\draw [line width=0.5pt] (-2,0)-- (2,0);
\draw [line width=0.5pt] (2,0)-- (1,-1);
\draw [line width=0.5pt] (1,-1)-- (-3,-1);
\draw [line width=0.5pt] (-3,-1)-- (-2,0);
\draw [rotate around={1.2383114065847913:(-0.4978893782460493,1.5108536116302471)},line width=0.5pt,fill=black,fill opacity=0.3] (-0.4978893782460493,1.5108536116302471) ellipse (0.5397691249101896cm and 0.1977822814176179cm);
\draw [rotate around={1.2383114065847909:(-0.5041514627319217,-0.5329809797708583)},line width=0.5pt,fill=black,fill opacity=0.3] (-0.5041514627319217,-0.5329809797708583) ellipse (0.5397691249101961cm and 0.19778228141761936cm);
\draw [samples=50,domain=-0.645:0.65,rotate around={0:(-0.5,0.5)},xshift=-0.5cm,yshift=0.5cm,line width=0.5pt] plot ({0.22238107315113698*(1+(\x)^2)/(1-(\x)^2)},{0.44782436099898426*2*(\x)/(1-(\x)^2)});
\draw [samples=50,domain=-0.645:0.65,rotate around={0:(-0.5,0.5)},xshift=-0.5cm,yshift=0.5cm,line width=0.5pt] plot ({0.22238107315113698*(-1-(\x)^2)/(1-(\x)^2)},{0.44782436099898426*(-2)*(\x)/(1-(\x)^2)});
\begin{scriptsize}
\draw (-1.8,-0.2) node {\large $E^{-}$};
\draw (-1.8,1.8) node {\large $E^{+}$};
\draw (-0.5,0.5) node {\large $K$};
\draw (-0.5,1.5) node {$B(0,R)$};
\end{scriptsize}
\end{tikzpicture}
\caption{Manifold $\mathcal{M}=\mathbb{R}^{n}\#\mathbb{R}^{n}$.}
\label{pic1}
\end{figure}

Let $n\geq 3$. Consider a manifold $\mathcal{M}=\mathbb{R}^{n}\#\mathbb{R}%
^{n}$ that is a connected sum of two copies of $\mathbb{R}^{n}$. This means
that $\mathcal{M}$ is a union of a compact part $K$ and two Euclidean ends $%
E^{\pm }=\mathbb{R}^{n}\backslash {B(0,R)}$ (see Fig. \ref{pic1}). 
This is a classical setting where the two-sided heat kernel estimate
\eqref{S2 twosided} fails, see for instance \cite{GrSa2009}.

Note that $V(x,R)\asymp {R}^{n}$ for all $x\in \mathcal{M}$, and the heat
kernel on $\mathcal{M}$ satisfies the upper bound (\ref{S2 upper}) (see, for
instance, \cite[Corollary 4.6, p.1946]{GrSa2009}). However,
the essential hypotheses, Hölder regularity \eqref{Holder} or pointwise 
gradient estimate \eqref{S2 gradient}, fail in the present setting 
(and so does the Poincaré inequality, \cite{GIS2021}).
Indeed, it was shown in  \cite[%
Proposition 6.1]{CCH2006} that, for large $t$, 
\begin{equation*}
\sup_{x,y\in \mathcal{M}}\left\vert \nabla _{y}h_{t}\left( x,y\right)
\right\vert \geq c\,t^{-n/2},
\end{equation*}%
which is incompatible with  (\ref{S2 gradient}). 

Let us verify that also the conclusion (\ref{S1 Main thm Linf}) of Theorem %
\ref{S1 Main thm} fails in this setting. For that consider on $E^{+}$ for any $t>1$
a point $x_{t}$ with $\left\vert x_{t}\right\vert =\sqrt{t}$ and with a
fixed direction $\omega =x_{t}/\sqrt{t}.$ It was proved in \cite[Proposition
6.1]{CCH2006} that, for any $x\in K,$ 
\begin{equation}
\lim_{t\rightarrow \infty }\,t^{\frac{n}{2}}h_{t}(x_{t},x)\,=\,\Phi (x)
\label{CCH}
\end{equation}%
where $\Phi \left( x\right) $ is a positive harmonic function on $\mathcal{M}
$ that  tends to a positive constant as $x\rightarrow \infty $ at the end $%
E^{+}$ and tends to $0$ as $x\rightarrow \infty $ at the end $E^{-}.$ Let us
fix two points $x_{1}$ and $x_{2}$ in $K$ such that $\Phi \left(
x_{1}\right) \neq \Phi \left( x_{2}\right) .$ Then (\ref{CCH}) implies that 
\begin{equation}
t^{\frac{n}{2}}\,\big(h_{t}(x_{t},x_{1})\,-\,h_{t}(x_{t},x_{2})\big)%
\,\longrightarrow \big(\Phi (x_{1})\,-\,\Phi (x_{2})\big)\qquad \textnormal{%
as}\;\;\;t\longrightarrow \infty   \label{x1x2}
\end{equation}%
while (\ref{S1 Main thm Linf}) would imply that%
\begin{equation*}
t^{\frac{n}{2}}\Vert h_{t}(\,\cdot\,,x_{1})\,-\,h_{t}(\,\cdot\,,x_{2})\Vert
_{L^{\infty}(\mathcal{M})}\longrightarrow 0
\qquad \textnormal{as}\;\;\;t\longrightarrow \infty\,,
\end{equation*}%
hence, 
\begin{equation*}
t^{\frac{n}{2}}\left(
h_{t}(x_{t},x_{1})\,-\,h_{t}(x_{t},x_{2})\right) \,\longrightarrow 0,
\end{equation*}%
which contradicts to (\ref{x1x2}).


\begin{acknowledgement}
The first author is funded by the Deutsche Forschungsgemeinschaft (DFG, German Research Foundation) - Project-ID 317210226 - SFB 1283. The second author is supported by the Hellenic Foundation for Research and Innovation, Project HFRI-FM17-1733. The last author acknowledges financial support from the Methusalem Programme \textit{Analysis and Partial Differential Equations (Grant number 01M01021)} during his postdoc stay at Ghent University.
\end{acknowledgement}
 
\input{references}

\end{document}

%% file: references.tex
%
%
%

%% file: GPZ2023.bbl
\begin{thebibliography}{99.}%
\bibitem{AnOs2003}
Anker, J.-Ph., Ostellari, P.: 
The heat kernel on noncompact symmetric spaces. 
In: Lie Groups And Symmetric Spaces, pp. 27--46. 
Amer. Math. Soc., Providence, RI (2003)

\bibitem{APZ2023}
Anker, J.-Ph., Papageorgiou, E., Zhang, H.-W.:
Asymptotic behavior of solutions to the heat equation on noncompact symmetric spaces. 
J. Funct. Anal. \textbf{284}, Paper No. 109828 (2023)

\bibitem{Bar1998}
Barlow, M.: 
Diffusions on fractals. In: Lectures On Probability Theory And Statistics (Saint-Flour, 1995), pp. 1--121.
Springer, Berlin (1998)

\bibitem{CCH2006}
Carron, G., Coulhon, T., Hassell, A.:
Riesz transform and {$L^p$}-cohomology for manifolds with Euclidean ends. 
Duke Math. J. \textbf{133}, 59-93 (2006)

\bibitem{ChHa2020}
Chen, X., Hassell, A.:
The heat kernel on asymptotically hyperbolic manifolds. 
Comm. Partial Differential Equations \textbf{45}, 1031-1071 (2020)


\bibitem{DaMa1988}
Davies, E., Mandouvalos, N.:
Heat kernel bounds on hyperbolic space and Kleinian groups. 
Proc. London Math. Soc. (3) \textbf{57}, 182-208 (1988)

\bibitem{Dav1990}
Davies, E.: 
Heat kernels and spectral theory. 
Cambridge University Press (1990)

\bibitem{Dun2004}
Dungey, N.:
Heat kernel estimates and Riesz transforms on some Riemannian covering manifolds.
Math. Z. \textbf{247}, 765-794 (2004)

\bibitem{DzPr2018}
Dziubański, J., Preisner, M.:
Hardy spaces for semigroups with Gaussian bounds. 
Ann. Mat. Pura Appl. (4). \textbf{197}, 
965-987 (2018)

\bibitem{FaSt1986}Fabes, E., Stroock, D.:
A new proof of Moser's parabolic Harnack inequality using the old ideas of Nash. 
Arch. Ration. Mech. Anal. \textbf{96}, 327-338 (1986)

\bibitem{Gri1991}
Grigor'yan, A.:
The heat equation on noncompact Riemannian manifolds. 
Mat. Sb. \textbf{182}, 55-87 (1991)

\bibitem{Gri1995}
Grigor'yan, A.: 
Upper bounds of derivatives of the heat kernel on an arbitrary complete manifold. 
J. Funct. Anal. \textbf{127}, 363-389 (1995)

\bibitem{Gri2009}Grigor'yan, A.:
Heat kernel and analysis on manifolds. 
AMS International Press (2009)

\bibitem{GrSa2009}
Grigor'yan, A., Saloff-Coste, L.:
Heat kernel on manifolds with ends. 
Ann. Inst. Fourier (Grenoble) \textbf{59}, 1917-1997 (2009)

\bibitem{GIS2021}
Grigor'yan, A., Ishiwata, S., Saloff-Coste, L.:
Geometric analysis on manifolds with ends. 
In: Analysis And Partial Differential Equations On Manifolds, Fractals And Graphs,
pp. 325-343, Adv. Anal. Geom. \textbf{3}, De Gruyter, Berlin (2021)

\bibitem{LiYa1986}
Li, P., Yau, S.:
On the parabolic kernel of the Schrödinger operator. 
Acta Math. \textbf{156}, 153-201 (1986)

\bibitem{Sal1992}
Saloff-Coste, L.:
A note on Poincaré, Sobolev, and Harnack inequalities.
Internat. Math. Res. Notices, 27-38 (1992)

\bibitem{Sal2002}
Saloff-Coste, L.:
Aspects of Sobolev-type inequalities.
Cambridge University Press (2002)

\bibitem{Sal2010}
Saloff-Coste, L.:
The heat kernel and its estimates.
In: Probabilistic Approach To Geometry, pp. 405-436.
Adv. Stud. Pure Math. \textbf{57}, Math. Soc. Japan, Tokyo (2010)

\bibitem{Str1983}
Strichartz, R.:
Analysis of the Laplacian on the complete Riemannian manifold.
J. Functional Analysi \textbf{52}, 48-79 (1983)

\bibitem{VSC1992}
Varopoulos, N. T., Saloff-Coste, L., Coulhon, T.:
Analysis and geometry on groups. 
Cambridge University Press (1992)

\bibitem{Vaz2018}
Vázquez, J.:
Asymptotic behaviour methods for the heat equation. Convergence to the Gaussian.
preprint, arXiv:1706.10034

\bibitem{Vaz2019}
Vázquez, J. Asymptotic behaviour for the heat equation in hyperbolic space.
To appear in Comm. Analysis and Geometry, arXiv:1811.09034

\end{thebibliography}
